\documentclass[11pt]{article}

\usepackage{amsmath}
\usepackage{amssymb,latexsym}
\usepackage{amsthm}
\usepackage{color}
\usepackage{graphicx}
\usepackage{mathabx}
\usepackage{hyperref}
\usepackage{fancyhdr}

\date{}
\newcommand{\cA}{{\mathcal A}}
\newcommand{\cB}{{\mathcal B}}

\newcommand{\cU}{{\mathcal U}}
\newcommand{\cP}{{\mathcal P}}
\newcommand{\cH}{{\mathcal H}}
\newcommand{\cD}{{\mathcal D}}
\newcommand{\cE}{{\mathcal E}}
\newcommand{\cL}{{\mathcal L}}
\newcommand{\cT}{{\mathcal T}}
\newcommand{\cC}{{\mathcal C}}
\newcommand{\GF}{\hbox{{\rm GF}}}
\newcommand{\so}{{\cal S}}

\renewcommand{\proof}{\noindent{\bf Proof.}\ }

\newcommand{\PG}{\mathrm{PG}}
\newcommand{\AG}{\mathrm{AG}}
\newcommand{\ep}{\epsilon}

\newtheorem{theorem}{Theorem}[section]
\newtheorem{lemma}[theorem]{Lemma}
\newtheorem{corollary}[theorem]{Corollary}
\newtheorem{remark}[theorem]{Remark}
\newtheorem{definition}[theorem]{Definition}
\newtheorem{proposition}[theorem]{Proposition}
\newtheorem{example}[theorem]{Example}
\newtheorem{conjecture}[theorem]{Conjecture}

\makeatother

\renewcommand{\qed}{\hfill \mbox{\raggedright \rule{0.1in}{0.1in}}}

\begin{document}

\title{On bisecants of R\'edei type blocking sets and applications}
\author{Bence Csajb\'{o}k\thanks{Research supported by the Hungarian National Foundation for Scientific Research, Grant No. K 81310 and by the Italian Ministry of Education, University and Research (PRIN 2012 project ``Strutture geometriche, combinatoria e loro applicazioni'')}}

\maketitle

\begin{abstract}
If $\mathcal{B}$ is a minimal blocking set of size less than $3(q+1)/2$ in $\mathrm{PG}(2,q)$, $q$ is a power of the prime $p$, then Sz\H{o}nyi's result states that each line meets $\mathcal{B}$ in $1 \pmod p$ points. It follows that $\mathcal{B}$ cannot have bisecants, i.e. lines meeting $\mathcal{B}$ in exactly two points. If $q>13$, then there is only one known minimal blocking set of size $3(q+1)/2$ in $\mathrm{PG}(2,q)$, the so called projective triangle. This blocking set is of R\'edei type and it has $3(q-1)/2$ bisecants, which have a very strict structure.
We use polynomial techniques to derive structural results on R\'edei type blocking sets from information on their bisecants.
We apply our results to point sets of $\mathrm{PG}(2,q)$ with few odd-secants.

In particular, we improve the lower bound of Balister, Bollob\'as, F\"uredi and Thompson on the number of odd-secants of a $(q+2)$-set in $\mathrm{PG}(2,q)$ and we answer a related open question of Vandendriessche. We prove structural results for semiovals and derive the non existence of semiovals of size $q+3$ when $p \neq 3$ and $q>5$. This extends a result of Blokhuis who classified semiovals of size $q+2$, and a result of Bartoli who classified semiovals of size $q+3$ when $q\leq 17$. In the $q$ even case we can say more applying a result of Sz\H{o}nyi and Weiner about the stability of sets of even type. We also obtain a new proof to a result of G\'acs and Weiner about $(q+t,t)$-arcs of type $(0,2,t)$ and to one part of a result of Ball, Blokhuis, Brouwer, Storme and Sz\H{o}nyi about functions over $\GF(q)$ determining less than $(q+3)/2$ directions.
\end{abstract}

\bigskip


\section{Introduction}

A \emph{blocking set} $\cB$ of $\PG(2,q)$, the Desarguesian projective plane of order $q$, is a point set meeting every line of the plane.
$\cB$ is called \emph{non-trivial} if it contains no line and \emph{minimal} if $\cB$ is minimal subject to set inclusion.
A point $P\in \cB$ is said to be \emph{essential} if $\cB \setminus \{P\}$ is not a blocking set.
For a point set $\so$ and a line $\ell$ we say that $\ell$ is a $k$-secant of $\so$ if $\ell$ meets $\so$ in $k$ points.
If $k=1$, $k=2$, or $k=3$, then we call $\ell$ a tangent to $\so$, a bisecant of $\so$, or a trisecant of $\so$, respectively.
We usually consider $\PG(2,q)$ as $\AG(2,q)$, the Desarguesian affine plane of order $q$, extended by the line at infinity, $\ell_{\infty}$. Throughout the paper $q$ will always denote a power of $p$, $p$ prime.
For the points of $\AG(2,q)$ we use cartesian coordinates. The infinite point (or \emph{direction}) of lines with slope $m$ will be denoted by $(m)$, the infinite point of vertical lines will be denoted by $(\infty)$. Let $\cU=\{(a_i,b_i)\}_{i=1}^q$ be a set of $q$ points of
$\AG(2,q)$. The set of \emph{directions determined by $\cU$} is $\cD_{\cU}:=\left\{\left(\frac{b_i-b_j}{a_i-a_j}\right) \colon i\neq j \right\}$.  It is easy to see that $\cB:=\cU \cup \cD_{\cU}$ is a blocking set of $\PG(2,q)$ with the property that there is a line, the line at infinity, which meets $\cB$ in exactly $|\cB|-q$ points. If $|\cD_{\cU}|\leq q$, then $\cB$ is minimal. Conversely, if $\cB$ is a minimal blocking set of size $q+N\leq 2q$ and there is a line meeting $\cB$ in $N$ points, then $\cB$ can be obtained from the above construction. Blocking sets of size $q+N\leq 2q$ with an $N$-secant are called blocking sets of \emph{R\'edei type}, the $N$-secants of the blocking set are called \emph{R\'edei lines}.
If the $q$-set $\cU$ does not determine every direction, then $\cU$ is affinely equivalent to the graph of a function $f$ from $\GF(q)$ to $\GF(q)$, i.e. $\cU=\{(x,f(x))\colon x\in \GF(q)\}$. Note that $f(x)-cx$ is a permutation polynomial if and only if $(c)$ is a direction not determined by the graph of $f$, see \cite{perpol} by Evans, Greene, Niederreiter. A blocking set is said to be \emph{small}, if its size is less than $q+(q+3)/2$. Small minimal R\'edei type blocking sets, or equivalently, functions determining less than $(q+3)/2$ directions,
have been characterized by Ball, Blokhuis, Brouwer, Storme, Sz\H{o}nyi and Ball, see \cite{BBBSSz, Ball}. From these results it follows that such blocking sets meet each line of the plane in $1 \pmod p$ points. This property holds for any small minimal blocking set, as it was proved by Sz\H{o}nyi in \cite{Szonyiblok}.

It follows from the above mentioned results that minimal blocking sets with bisecants cannot be small. If $q$ is odd, then the smallest known non-small minimal R\'edei type blocking set is the following set of
$q+(q+3)/2$ points (up to projective equivalence):
\[\cB:=\{(0:1:a),\,(1:0:a),\,(-a:1:0) \colon \text{$a$ a square in $\GF(q)$}\} \cup \{(0:0:1)\}.\]
In the book of Hirschfeld {\cite[Lemma 13.6 (i)]{jwph}} this example is called the \emph{projective triangle}.
$\cB$ has three R\'edei lines and has the following properties. Through each point of $\cB$ there passes a bisecant of $\cB$. If $\cH \subset \cB$ is a set of collinear points such that there passes a unique bisecant of $\cB$ through each point of $\cH$ and there is a R\'edei line $\ell$ disjoint from $\cH$, then the bisecants through the points of $\cH$ are contained in a pencil. In Theorem \ref{main1} we show that this property holds for any R\'edei type blocking set. In fact, we prove the following stronger result.
If $R_1$ and $R_2$ are points of $\cB \setminus \ell$, such that for $i=1,2$ there is a unique bisecant of $\cB$ through $R_i$ and there is
a point $T\in \ell$, such that $TR_1$ and $TR_2$ meet $\cB$ in at least four points, then for each $M\in \ell$ the lines $R_1M$ and $R_2M$ meet $\cB$ in the same number of points.
The essential part of our proof is algebraic, it is based on polynomials over $\GF(q)$.
We apply our results to point sets of $\PG(2,q)$ with few odd-secants, which we detail in the next paragraphs.

A semioval $\so$ of a finite projective plane is a point set with the property that at each point of $\so$ there passes exactly one tangent to $\so$. For a survey on semiovals see \cite{survey} by Kiss.
In $\PG(2,q)$ Blokhuis characterized semiovals of size $q-1+a$, $a>2$,  meeting each line in 0,1,2, or $a$ points. He also proved that there is no semioval of size $q+2$ in $\PG(2,q)$, $q>7$, see \cite{seminuclear} and \cite{originsemi}, where the term \emph{seminuclear set}  was used for semiovals of size $q+2$. For another characterization of semiovals with special intersection pattern with respect to lines see \cite{gacssemi} by G\'acs. We refine Blokhuis' characterization to obtain new structural results about semiovals of size $q-1+a$ containing
$a$ collinear points. As an application, we prove the non-existence of semiovals of size $q+3$ in $\PG(2,q)$, $5<q$ odd when $p\neq 3$. For $q\leq 17$ this was also proved by Bartoli in \cite{db}. When $q$ is small, then the spectrum of the sizes of semiovals in $\PG(2,q)$ is known, see \cite{lisonek} by Lisonek for $q\leq 7$ and \cite{smalls} by Kiss, Marcugini and Pambianco for $q=9$.
When $q$ is even, then a stronger result follows from \cite[Theorem 5.3]{setofeven} by Sz\H{o}nyi and Weiner on the stability of sets of even type.

In the recent article \cite{symmdiff} by Balister, Bollob\'{a}s, F\"{u}redi and Thompson, the minimum number of odd-secants of an $n$-set in
$\PG(2,q)$, $q$ odd, was investigated. They studied in detail the case of $n=q+2$. In our last section we improve their lower bound and
we answer a related open question of Vandendriessche from \cite{PV}.

Our Theorem \ref{main0} yields a new proof to \cite[Theorem 2.5]{gw} by G\'acs and Weiner about $(q+t,t)$-arcs of type $(0,2,t)$.
In Section \ref{dirp} we explain some connections between Theorem \ref{main0} and the direction problem.

\section{Bisecants of R\'edei type blocking sets}

\begin{lemma}
\label{poly}
Let $\cU$ be a set of $q$ points in $\AG(2,q)$ and denote by $\cD_{\cU}$ the set of directions determined by $\cU$.
Take a point $R=(a_0,b_0)\in \cU$ and denote the remaining $q-1$ points of $\cU$ by $(a_i,b_i)$ for $i=1,2,\ldots,q-1$.
Consider the following polynomial:
\begin{equation}
\label{pp}
f(Y):=\prod_{i=1}^{q-1}((a_i-a_0)Y-(b_i-b_0))\in \GF(q)[Y].
\end{equation}
For $m\in \GF(q)$ the following holds.
\begin{enumerate}
	\item The line through $R$ with direction $m$ meets $\cU$ in $k_m$ points if and only if $m$ is a $(k_m-1)$-fold root of $f(Y)$.
	\item If $(m)\notin \cD_{\cU}$, then $f(m)=-1$.
	\item If $(\infty)\notin \cD_{\cU}$, then the coefficient of $Y^{q-1}$ in $f$ is $-1$.
\end{enumerate}
\end{lemma}
\proof
We have $(a_j-a_0)m-(b_j-b_0)=0$ for some $j\in \{1,2,\ldots,q-1\}$ if and only if $(m)$, $R$ and $(a_j,b_j)$ are collinear.
This proves part 1. To prove part 2, note that $(a_j-a_0)m-(b_j-b_0)=(a_k-a_0)m-(b_k-b_0)$ for some
$j,k\in\{1,2,\ldots,q-1\}$, $j\neq k$, if and only if $(a_j-a_k)m-(b_j-b_k)=0$, i.e. if and only if $(a_j,b_j)$, $(a_k,b_k)$ and $(m)$ are collinear. If $(m)\notin \cD_{\cU}$, then this cannot be and hence $\{(a_i-a_0)m-(b_i-b_0) \colon i=1,2,\ldots,q-1\}$ is the set of non-zero elements of $\GF(q)$. It follows that in this case $f(m)=-1$. If $(\infty)\notin \cD_{\cU}$, then $\{a_i-a_0 \colon i=1,2,\ldots,q-1\}$ is the set of non-zero elements of $\GF(q)$, and hence $\prod_{i=1}^{q-1}(a_i-a_0)=-1$.
\qed

\begin{remark}
For a set of affine points $\cU=\{(a_i,b_i)\}_{i=0}^{k}$ the R\'edei polynomial of $\cU$ is
$\prod_{i=0}^k(X+a_iY-b_i)=\sum_{j=0}^{k+1} h_j(Y)X^{k+1-j}\in \GF(q)[X,Y]$,
where $h_j(Y)\in\GF(q)[Y]$ is a polynomial of degree at most $j$.
Now suppose that $\cU$ is a $q$-set and $(a_0,b_0)=(0,0)$. Then
$h_{q-1}(Y)=\sum_{j=0}^{q-1} \prod_{i\neq j} (a_iY-b_i)=\prod_{i=1}^{q-1}(a_iY-b_i)$ is the polynomial
associated to the affine $q$-set $\cU$ as in Lemma \ref{poly}.
This polynomial also appears in Section 4 of Ball's paper \cite{Ball}.
\end{remark}

\begin{theorem}
\label{main0}
Let $\cB$ be a blocking set of R\'edei type in $\PG(2,q)$, with R\'edei line $\ell$.
\begin{enumerate}
	\item If there is a point in $\cB\setminus \ell$ which is not incident with any bisecant of $\cB$, then $\cB$ is minimal and $|\ell\cap \cB|\equiv 1 \pmod p$.
	\item If $R,R'\in \cB\setminus \ell$ such that $R$ and $R'$ are not incident with any bisecant of $\cB$, then $|RM \cap \cB|=|R'M\cap \cB|$ for each $M\in \ell$.
\end{enumerate}
\end{theorem}
\proof
It is easy to see that if there is a point $R\in \cB\setminus \ell$, such that there is no bisecant of $\cB$ through $R$, then
$|\cB \cap \ell|\leq q-1$. First we show that $\cB$ is minimal.
As $\cB$ is of R\'edei type, the points of $\cB \setminus \ell$ are essential in $\cB$.
Take a point $D\in \cB \cap \ell$. As there is no bisecant through $R$, it follows that $DR$ meets $\cB$ in at least three points
and hence there is a tangent to $\cB$ at $D$, i.e. $D$ is essential in $\cB$.

We may assume that $\ell=\ell_{\infty}$ and $(\infty)\notin \cB$. Let $R=(a_0,b_0)$ be a point of $\cB\setminus \ell$ which is not incident with any bisecant of $\cB$ and let
$\cU=\cB\setminus \ell_{\infty}=\{(a_i,b_i)\}_{i=0}^{q-1}$. Consider the polynomial $f(Y)=\prod_{i=1}^{q-1}((a_i-a_0)Y-(b_i-b_0))$ introduced in
\eqref{pp}. Let $m\in \GF(q)$. According to Lemma \ref{poly} we have the following.
\begin{itemize}
\item If $(m)\in \cB$, then $f(m)=0$,	
\item if $(m)\notin \cB$, then $f(m)=-1$,
\item the coefficient of $Y^{q-1}$ in $f$ is $-1$.
\end{itemize}
Now let $\ell_{\infty} \setminus (\cB\cup \{(\infty)\})=\{(m_1),\,(m_2),\ldots,\,(m_k)\}$ and consider the polynomial
\[g(Y):=\sum_{i=1}^k(Y-m_i)^{q-1}-k.\]
For $m\in \GF(q)$ we have $g(m)=f(m)$. As both polynomials have degree at most $q-1$, it follows that
$g(Y)=f(Y)$. The coefficient of $Y^{q-1}$ is $k$ in $g$ and hence $p \divides k+1$.
As $k+1 =q+1-|\cB\cap \ell_{\infty}|$, part 1 follows.

For $(m)\notin \cB$ the line through any point of $\cU$ with slope $m$ meets $\cB$ in $1$ point.
For $(m)\in \cB$ the line through $R$ with slope $m$ meets $\cB$ in $k_m+2$ points if and only if $m$ is a
$k_m$-fold root of $f(Y)$. As $f(Y)=g(Y)$, and the coefficients of $g(Y)$ depend only on the points
of $\cB\cap \ell_{\infty}$, it follows that $k_m$ does not depend on the initial choice of the point $R$,
as long as the chosen point is not incident with any bisecant of $\cB$. This proves part 2.
\qed

\begin{theorem}
\label{main1}
Let $\cB$ be a blocking set of R\'edei type in $\PG(2,q)$, with R\'edei line $\ell$.
\begin{enumerate}
	\item If there is a point in $\cB\setminus \ell$ contained in a unique bisecant of $\cB$, then $|\cB \cap \ell|\not\equiv 1 \pmod p$.
	\item If $R_1,R_2\in \cB \setminus \ell$, each of them is contained in a unique bisecant of $\cB$ and there is a point $T\in \ell$
	such that $R_1T$ and $R_2T$ both meet $\cB$ in at least four points, then for each $M\in \ell$ we have $|MR_1 \cap \cB|=|MR_2 \cap \cB|$.
	\item If $R_1,R_2\in \cB \setminus \ell$, each of them is contained in a unique bisecant of $\cB$ and the common point of these bisecants is on the line $\ell$, then for each $M\in \ell$ we have $|MR_1 \cap \cB|=|MR_2 \cap \cB|$.
\end{enumerate}
\end{theorem}
\proof
Let $R$ be a point of $\cB\setminus \ell$ contained in a unique bisecant $r$ of $\cB$.
First suppose $|\cB \cap \ell|=q$. Then part 1 is trivial and there is no line through $R$ meeting $\cB$ in at least 4 points, since otherwise we would get more than one bisecants through $R$. Suppose that $R'$ is another point of $\cB\setminus \ell$ contained in a unique bisecant $r'$ of $\cB$ and $r\cap r' \in \ell$. Let $\{Q\}=\ell \setminus \cB$. Then $RQ$ and $R'Q$ are tangents to $\cB$ and $|MR \cap \cB|=|MR'\cap \cB|=3$ for each $M\in (\ell \cap \cB) \setminus \{r\cap r'\}$.
From now on, we assume $k:=q-|\cB \cap \ell|\geq 1$.

First we prove the theorem when $\cB$ is minimal.
We may assume $\ell=\ell_{\infty}$ and
$\ell_{\infty} \setminus \cB = \{(\infty),(m_1),\ldots,(m_k)\}$.

As in the proof of Theorem \ref{main0}, let $\cU=\cB\setminus \ell_{\infty}=\{(a_i,b_i)\}_{i=0}^{q-1}$ and define $f(Y)$ as in \eqref{pp}.
Take $m\in \GF(q)$ and let $t$ be the slope of the unique bisecant through $R$. From Lemma \ref{poly} we obtain the following.
\[
f(m)=\left\{ \begin{array}{ll}
-1 & \textrm{if $(m)\notin \cB$}, \\
0  & \textrm{if $(m)\in \cB \setminus \{(t)\}$}, \\
f(t)\neq 0 & \textrm{if $m=t$}.
\end{array} \right.
\]
Consider the polynomial
\begin{equation}
\label{eq3}
g(Y):=f(t)+|\cB\cap \ell_{\infty}|+\sum_{i=1}^{k}(Y-m_i)^{q-1}-f(t)(Y-t)^{q-1}.
\end{equation}
For $m\in \GF(q)$ we have $g(m)=f(m)$. As both polynomials have degree at most $q-1$, it follows that
$g(Y)=f(Y)$. The coefficient of $Y^{q-1}$ is $-|\cB\cap \ell_{\infty}|-f(t)$ in $g$ and $-1$ in
$f$. It follows that $p \divides |\cB\cap \ell_{\infty}|+f(t)-1$ and hence $f(t)\equiv 1-|\cB\cap \ell_{\infty}|\equiv k+1 \pmod p$.
If $|\cB\cap \ell_{\infty}| \equiv 1 \pmod p$, then $f(t)=0$, a contradiction. This proves part 1.

Now consider
\[\partial_Y g(Y)=-\sum_{i=1}^{k}(Y-m_i)^{q-2}+(k+1)(Y-t)^{q-2},\]
and
\[w(Y):=(Y-t)\prod_{i=1}^k (Y-m_i)\partial_Yg(Y)=\]
\[-\sum_{i=1}^{k}(Y-m_i)^{q-1}(Y-t)\prod_{j\neq i}(Y-m_j)+(k+1)(Y-t)^{q-1}\prod_{j=1}^k(Y-m_j).\]
If $(m)\in \cB \setminus \{(t)\}$, then
\[w(m)=-\sum_{i=1}^{k}(m-t)\prod_{j\neq i}(m-m_j)+(k+1)\prod_{j=1}^k(m-m_j).\]
Suppose that the line through $R$ with direction $m$ meets $\cB$ in at least four points. Then $m$ is a multiple root of $f(Y)$ and hence
it is also a root of $w(Y)$. It follows that $m$ is a root of
\begin{equation}
\tilde{w}(Y):=-(Y-t)\sum_{i=1}^{k}\prod_{j\neq i}(Y-m_j)+(k+1)\prod_{j=1}^k(Y-m_j).
\end{equation}
Note that $\sum_{i=1}^{k}\prod_{j\neq i}(m-m_j)=0$ and $\tilde{w}(m)=0$ together would imply $(k+1)\prod_{j=1}^k(m-m_j)=0$, which cannot be since $(m)\notin \{(m_1),\ldots,(m_k)\}$ and $p\notdivides k+1$. It follows that $t$ can be expressed from $m$ and $m_1,\ldots,m_k$ in the following way:
\begin{equation}
\label{eqq}
t=m-\frac{(k+1)\prod_{j=1}^k(m-m_j)}{\sum_{i=1}^{k}\prod_{j\neq i}(m-m_j)}.
\end{equation}

Now let $R_1$ and $R_2$ be two points as in part 2 and let $T=(m)$. It follows from \eqref{eqq} that the bisecants through these points have the same slope. Then, according to \eqref{eq3}, $f(Y)=g(Y)$ does not depend on the choice of $R_i$, for $i=1,2$. The assertion follows from Lemma \ref{poly} part 1.

If $R_1$ and $R_2$ are two points as in part 3, then the bisecants through these points have the same slope. It follows that
$f(Y)=g(Y)$ does not depend on the choice of $R_i$, for $i=1,2$. As above, the assertion follows from Lemma \ref{poly} part 1.

Now suppose that $\cB$ is not minimal and $R_1\in \cB\setminus \ell$ is contained in a unique bisecant of $\cB$.
As $\cB$ is a blocking set of R\'edei type, the points of $\cB \setminus \ell$ are essential in $\cB$.
Let $C\in \cB \cap \ell$ such that $\cB':=\cB \setminus \{C\}$ is a blocking set.
In this case for each $P\in \cB \setminus \ell$ the line $PC$ is a bisecant of $\cB$ and $R_1C$ is the unique bisecant of $\cB$ through $R_1$.
It follows that there is no bisecant of $\cB'$ through $R_1$. Then Theorem \ref{main0} yields that $|\ell \cap \cB'|\equiv 1 \pmod p$.
As $|\ell \cap \cB|=|\ell \cap \cB'|+1$, we proved part 1.

If $R_2$ is another point of $\cB\setminus \ell$ such that $R_2$ is contained in a unique bisecant of $\cB$, then there is no bisecant of $\cB'$ through $R_2$ and hence parts 2 and 3 follow from Theorem \ref{main0} part 2.
\qed

\section{Connections with the direction problem}
\label{dirp}

Let $\cB$ be a blocking set in $\PG(2,q)$. We recall $q=p^h$, $p$ prime. The \emph{exponent} of $\cB$ is the maximal integer $0 \leq e \leq h$
such that each line meets $\cB$ in $1 \pmod {p^e}$ points.
We recall the following two results about the exponent.

\begin{theorem}[Sz\H{o}nyi \cite{Szonyiblok}]
\label{modpe}
Let $\cB$ be a small minimal blocking set in $\mathrm{PG}(2,q)$. Then $\cB$ has positive exponent.
\end{theorem}

\begin{theorem}[Sziklai \cite{linconj}]
\label{Sziklai}
Let $\cB$ be a small minimal blocking set in $\mathrm{PG}(2,q)$. Then the exponent of $\cB$ divides $h$.
\end{theorem}

\begin{proposition}
\label{cor0}
Let $\cB$ be a blocking set of R\'edei type in $\PG(2,q)$, with R\'edei line $\ell$ .
Suppose that $\cB$ does not have bisecants. Then $\cB$ has positive exponent and for each point $M\in \ell\cap \cB$ the lines through $M$ different from $\ell$  meet $\cB$ in 1 or in $p^t+1$ points, where $t$ is a positive integer depending only on the choice of $M$.
\end{proposition}
\proof
Theorem \ref{main0} part 1 yields that $\ell$ meets $\cB$ in $1 \pmod p$ points.
Lines meeting $\ell$ not in $\cB$ are tangents to $\cB$. For any $M\in \ell \cap \cB$
Theorem \ref{main0} part 2 yields that $MR$ meets $\cB \setminus \ell$ in the same
number of points for each $R\in \cB \setminus \ell$. Denote this number by $k$.
Then $k$ divides $|\cB \setminus \ell|=q$. As $\cB$ does not have bisecants, it follows that $k>1$ and hence $k=p^t$ for some positive integer $t$.
\qed

The following result is a consequence of the lower bound on the size of an affine blocking set due to Brouwer and Schrijver \cite{BS} and Jamison \cite{J}.

\begin{theorem}[Blokhuis and Brouwer {\cite[pg. 133]{BB}}]
\label{tg}
If $\cB$ is a minimal blocking set of size $q+N$, then there are at least $q+1-N$ tangents to $\cB$ at each point of $\cB$.
\end{theorem}

\begin{theorem}
\label{func}
Let $f$ be a function from $\GF(q)$ to $\GF(q)$ and let $N$ be the number of directions determined by $f$. If any line with a direction determined by $f$ that is incident with a point of the graph of $f$ is incident with at least two points of the graph of $f$, then each line meets the graph of $f$ in $p^t$ points for some integer $t$ and
\[q/s + 1 \leq N \leq (q - 1)/(s - 1),\]
where $s=\min \{p^t \colon \, \text{there is line meeting the graph of $f$ in $p^t > 1$ points} \}$.
\end{theorem}
\proof
If $\cU$ denotes the graph of $f$, then $\cB:=\cU \cup \cD_{\cU}$ is a blocking set of R\'edei type without bisecants.
Proposition \ref{cor0} yields that each line meets $\cU$ in $p^t$ points for some integer $t$, with $t=0$ only for lines with direction not in $\cD_{\cU}$.
Take a point $R\in \cU$ and let $\cD_{\cU}=\{D_1,D_2,\ldots, D_N\}$.
Then $|D_iR \cap \cB|\geq s+1$ yields $|\cB|=q+N\geq Ns+1$ and hence $(q-1)/(s-1)\geq N$.
Take a line $m$ meeting $\cU$ in $s$ points and let $M=m\cap \ell_{\infty}$.
According to Proposition \ref{cor0} the lines through $M$ meet $\cU$ in 0 or in $s$ points. Theorem \ref{tg}
yields that the number of lines through $M$ that meet $\cU$ is at most $N-1$. It follows that
$(N-1)s\geq q$ and hence $N\geq q/s+1$.
\qed
\medskip

Applying Theorems \ref{func} and \ref{modpe} we can give a new proof to the following result.

\begin{theorem}[part of Ball et al. \cite{BBBSSz} and Ball \cite{Ball}]
\label{part}
Let $f$ be a function from $\GF(q)$ to $\GF(q)$ and let $N$ be the number of directions determined by $f$.
Let $s=p^e$ be maximal such that any line with a direction determined by $f$ that is incident with a point of the graph of $f$ is incident with a multiple of $s$ points of the graph of $f$. Then one of the following holds.
\begin{enumerate}
	\item $s=1$ and $(q+3)/2 \leq N \leq q+1$,
	\item $q/s + 1 \leq N \leq (q - 1)/(s - 1)$,
	\item $s=q$ and $N=1$.
\end{enumerate}
\end{theorem}
\proof
The point set $\cB:=\cU \cup \cD_{\cU}$ is a minimal blocking set of R\'edei type.
If $s=1$, then $\cB$ cannot be small because of Sz\H{o}nyi's Theorem \ref{modpe} and hence $N\geq (q+3)/2$.
If $s>1$, then the bounds on $N$ follow from Theorem \ref{func}.
\qed
\medskip

In \cite{BBBSSz} and \cite{Ball} it was also proved that for $s>2$ the graph of $f$ is $\GF(s)$-linear and
that $\GF(s)$ is a subfield of $\GF(q)$. Note that Theorem \ref{Sziklai} by Sziklai generalizes the latter result.

\section{Small semiovals}

An \emph{oval} of a projective plane of order $q$ is a set of $q+1$ points such that no three of them are collinear.
It is easy to see that ovals are semiovals. The smallest known \emph{non-oval semioval}, i.e. semioval which is not an oval, is due to Blokhuis.

\begin{example}[Blokhuis \cite{seminuclear}]
\label{Blokhuis}
Let $\so$ be the following point set in $\PG(2,q)$, $3<q$ odd,
$\so=\{(0:1:s),\,(s:0:1),\,(1:s:0)\colon \text{ $-s$ is not a square}\}.$
Then $\so$ is a semioval of size $3(q-1)/2$.
\end{example}

\begin{conjecture}[Kiss et al. {\cite[Conjecture 11]{smalls}}]
If a semioval in $\PG(2,q)$, $q>7$, has less than $3(q-1)/2$ points, then it has exactly $q+1$ points and it is an oval.
\end{conjecture}

Let $\so$ be a semioval and $\ell$ a line meeting $\so$ in at least two points.
Take a point $P\in \so \cap \ell$. As there is a unique tangent to $\so$ at $P$, it follows that $|\so \setminus \ell|\geq q-1$, and hence
$|\so|\geq |\so \cap \ell|+q-1\geq q+1$. It is convenient to denote the size of $\so$ by $q-1+a$, where $a\geq 2$ holds automatically.
Then each line meets $\so$ in at most $a$ points.

\begin{theorem}[Blokhuis \cite{seminuclear}]
\label{Blokhuis2}
Let $\so$ be a semioval of size $q-1+a$, $a>2$, in $\PG(2,q)$ and suppose that each line meets $\so$ in $0$, $1$, $2$, or in $a$ points.
Then $\so$ is the symmetric difference of two lines with one further point removed from both lines, or $\so$ is projectively equivalent to Example \ref{Blokhuis}.
\end{theorem}

If $\so$ is a semioval of size $q+2$, then each line meets $\so$ in at most three points, thus Theorem \ref{Blokhuis2} yields the following.

\begin{theorem}[Blokhuis \cite{seminuclear}]
\label{Blokhuis3}
Let $\so$ be a semioval of size $q+2$ in $\PG(2,q)$. Then $\so$ is the symmetric difference of two lines with one further point removed
from both lines in $\PG(2,4)$, or $\so$ is projectively equivalent to Example \ref{Blokhuis} in $\PG(2,7)$.
\end{theorem}

We also recall the following well-known result by Blokhuis which will be applied several times. For another proof and possible generalizations see \cite[Remark 7]{dirszonyi} by Sz\H{o}nyi, or \cite[Corollary 3.6]{semiar2} by Csajb\'{o}k, H\'eger and Kiss.

\begin{proposition}[Blokhuis {\cite[Proposition 2]{seminuclear}}]
\label{nucleus}
Let $\so$ be a point set of $\PG(2,q)$, $q>2$, of size $q-1+a$, $a\geq 2$, with an $a$-secant $\ell$.
If there is a unique tangent to $\so$ at each point of $\ell\cap \so$, then these tangents are contained in a pencil.
The carrier of this pencil is called the \emph{nucleus} of $\ell$ and it is denoted by $N_{\ell}$.
For the sake of simplicity, the nucleus of a line $\ell_i$ will be denoted by $N_i$.
\end{proposition}

If $\cA$ and $\cB$ are two point sets, then $\cA \Delta \cB$ denotes their symmetric difference, that is $(\cA\setminus \cB)\cup (\cB\setminus \cA)$.

\begin{example}[Csajb\'ok, H\'eger and Kiss {\cite[Example 2.12]{semiar2}}]
\label{altered}
Let $\cB'$ be a blocking set of R\'{e}dei type in $\PG(2,q)$, with R\'edei line $\ell$.
Suppose that there is a point $P\in \cB' \setminus \ell$ such that the bisecants of $\cB'$ pass through $P$ and
there is no trisecant of $\cB'$ through $P$.
For example, if $\cB'$ has exponent $e$ and $p^e\geq 3$ (cf. Section \ref{dirp}), then $\cB'$ has no bisecants or trisecants and hence one can choose any point $P\in \cB' \setminus \ell$.
Take a point $W\in \ell \setminus \cB'$ and let $\so = (\ell \Delta \cB') \setminus \{W,P\}$.
Then $\so$ is a semioval of size $q-1+a$, where $a=|\ell \cap \so|$.
\end{example}

\begin{remark}
The blocking set $\cB'$ in Example \ref{altered} is necessarily minimal. To see this consider any point $R\in \cB' \setminus (\ell\cup \{P\})$.
As the bisecants of $\cB'$ pass through $P$, it follows that there is no bisecant of $\cB'$ through $R$ and hence Theorem \ref{main0} part 1  yields that $\cB'$ is minimal. \qed
\end{remark}

\begin{lemma}
\label{semioval}
Let $\so$ be a semioval of size $q-1+a$ in $\PG(2,q)$ and suppose that there is a line $\ell$ which is an $a$-secant of $\so$.
Denote the set of tangents through the points of $\so \setminus \ell$ by $\cL$ and let $\cB=\{N_{\ell}\}\cup (\so \Delta \ell)$. Then one of the following holds.
\begin{enumerate}
	\item $\so$ is an oval.
	\item $\cL$ is contained in a pencil with carrier $C$.
	Then $C\in\ell$ and $\cB':=\cB \setminus \{C\}$ is a blocking set of R\'edei type with R\'edei line $\ell$.
	In this case $\so$ can be obtained from $\cB'$ as in Example \ref{altered} with $P=N_{\ell}$ and $W=C$.
	\item $\cL$ is not contained in a pencil. Then $\cB$ is a minimal blocking set of R\'edei type with R\'edei line $\ell$ and
	\begin{enumerate}
		\item $p \notdivides a$,
		\item for any $R\in \so \setminus \ell$ the line $RN_{\ell}$ is not a tangent to $\so$,
		\item if $R_1,R_2\in \so\setminus \ell$ and there is a point $T\in \ell$ such that $R_iT$ meets $\so \cup \{N_{\ell}\}$ in at least three points for $i=1,2$, then for each $M\in \ell$ we have $|R_1M \cap (\so \cup \{N_{\ell}\})|=|R_2M\cap (\so\cup \{N_{\ell}\})|$,
		\item if $R_1,R_2\in \so\setminus \ell$ and the tangents to $\so$ at these two points meet each other on the line $\ell$, then for each $M\in \ell$ we have $|R_1M \cap (\so \cup \{N_{\ell}\})|=|R_2M\cap (\so\cup \{N_{\ell}\})|$.
	\end{enumerate}
	\end{enumerate}
\end{lemma}
\proof
First we show that $\cB$ is a blocking set of R\'edei type.
Take a point $R\in \so \setminus \ell$. As there is a tangent to $\so$ at $R$ it follows that $\ell$ meets $\so$ in at most $q$ points and
hence $\ell$ is blocked by $\cB$. Lines meeting $\ell$ not in $\so$ are blocked by $\cB$ since $\ell \setminus \so \subset \cB$.
If a line $m$ meets $\ell$ in $\so$, then either $m$ is a tangent to $\so$ and hence $N_{\ell}\in m$, or $m$ is not a tangent to $\so$ and hence there is a point of $\so \setminus \ell$ contained in $m$. As $\{N_{\ell}\} \cup (\so \setminus \ell) \subset \cB$, it follows that
$m$ is blocked by $\cB$ and hence $\cB$ is a blocking set. The line $\ell$ meets $\cB$ in $|\cB|-q$ points, thus $\cB$ is of R\'edei type and $\ell$ is a R\'edei line of $\cB$.

If $a=2$, then $\so$ is an oval. From now on we assume $a\geq 3$.
First suppose that $\cL$ is contained in a pencil with carrier $C$.
If $C \notin \ell$, then $|\cL|\leq q+1-a$, but $|\cL|=|\so \setminus \ell|=q-1$. It follows that $C\in \ell$.

Let $\cB'=\cB \setminus \{C\}$. In this paragraph we prove that $\cB'$ is a blocking set. It is enough to show that the lines through $C$ are blocked by $\cB'$. This trivially holds for the $q-1$ lines in $\cL$. First we show that $\cB'$ blocks $\ell$ too. Suppose to the contrary that $\ell \setminus (\so \cup \{C\})=\emptyset$ and hence $a=q$. As $a\geq 3$, we have $q\geq 3$ and hence there are at least two points in $\so \setminus \ell$. Take $R,Q\in \so \setminus \ell$ and let $M=RQ \cap \ell$. Since $M\neq C$, we have $M\in \so$. Then there are at least two tangents to $\so$ incident with $M$ and this contradiction shows that $\ell$ is blocked by $\cB'$. Now we show $CN_{\ell} \notin \cL$. Suppose to the contrary that $CN_{\ell}$ is a tangent to $\so$ at some $V\in \so \setminus \ell$. Then $VC$ is a trisecant of $\cB$.
If there were a bisecant $v$ of $\cB$ through $V$, then, by the construction of $\cB$, $v$ would be a tangent to $\so$ at $V$. This cannot be since the unique tangent to $\so$ at $V$ is $VC$, which is a trisecant of $\cB$ and hence $v\neq VC$.
For any $V'\in \so \setminus (\ell \cup \{V\})$, there is a unique bisecant of $\cB$ through $V'$, namely $V'C$.
We have shown that there is a point in $\cB\setminus \ell$ not incident with any bisecant of $\cB$ and there are points in $\cB \setminus \ell$ incident with a  unique bisecant of $\cB$. This cannot be because of Theorem \ref{main0} part 1 and Theorem \ref{main1} part 1.
It follows that $CN_{\ell}$ is not a tangent to $\so$.
As $CN_\ell$ is blocked by $\cB'$ and the other $q$ lines through $C$, $\ell$ and the lines of $\cL$, are also blocked, it follows that $\cB'$ is a blocking set. It is easy to see that $\ell$ is a R\'edei line of $\cB'$.

We show that there is no bisecant of $\cB'$ through the points of $\so \setminus \ell$.
Take a point $R\in \so\setminus \ell$ and suppose to the contrary that there is a bisecant $b$ of $\cB'$ through $R$.
Then, by the construction of $\cB'$, the line $b$ is a tangent to $\so$ at $R$. This is a contradiction since $b\neq RC$. It follows that if $\cB'$ has bisecants, then they pass through $N_{\ell}$.
If there were a trisecant $t$ of $\cB'$ through $N_{\ell}$, then let $V=t \cap \so$. It follows that $t$ is a tangent to $\so$ at $V$.
But we have already seen that there is no line of $\cL$ incident with $N_{\ell}$. This finishes the proof of part 2.

Now suppose that $\so$ is as in part 3.
If $\cB$ were not minimal, then the line set $\cL$ would be contained in a pencil with carrier on $\ell$, a contradiction.
Take a point $R\in \so\setminus\ell$.
If $RN_{\ell}$ is the tangent to $\so$ at $R$, then there is no bisecant of $\cB$ through $R$, thus $p\divides a$ (cf. Theorem \ref{main0} part 1). If $RN_{\ell}$ is not the tangent to $\so$ at $R$, then there is a unique bisecant of $\cB$ through $R$ (the tangent to $\so$ at $R$), thus $p\notdivides a$ (cf. Theorem \ref{main1} part 1). It follows that if any of the lines of $\cL$ is incident with $N_{\ell}$, or if $p\divides a$, then the whole line set $\cL$ is contained in the pencil with carrier $N_{\ell}$, a contradiction. This proves parts (a) and (b). Parts (c) and (d) follow from Theorem \ref{main1} parts 2 and 3, respectively.
\qed

\begin{remark}
\label{rem}
The properties (a)-(d) in part 3 of Lemma \ref{semioval} also hold when $\so$ is as in Example \ref{altered}.
From the properties of the point $P$ in Example \ref{altered} it follows that for $R\in \so \setminus \ell$ the line $RP$ is not a
tangent to $\so$ and this proves (b). As for any two points $R_1, R_2 \in \so \setminus \ell$ there is no bisecant of $\cB'$ incident with $R_1$ or $R_2$, properties (a), (c) and (d) follow from Theorem \ref{main0}. \qed
\end{remark}

\begin{theorem}
\label{semioval2}
Let $\so$ be a semioval of size $q-1+a$, $a>2$, which admits an $a$-secant $\ell$, and let $m\neq \ell$ be a $k$-secant of $\so$.
\begin{enumerate}
		\item For each $R\in \so \setminus \ell$, the line $RN_{\ell}$ is not a tangent to $\so$.
		\item If $k\geq 3$, then the tangents to $\so$ at the points of $m$ are contained in a pencil with carrier on $\ell$.
		\item If $k>(a-1)/2$, then $k=a$ and $N_{\ell}\in m$, or $k=\left\lceil a/2\right\rceil$ and $N_{\ell}\notin m$.
\end{enumerate}
\end{theorem}
\proof
Part 1 follows from Lemma \ref{semioval} part 3 (b), and part 2 follows from Lemma \ref{semioval} part (c) with $T=m\cap \ell$.

To prove part 3 first suppose $k>(a+1)/2$ and $N_{\ell}\notin m$. Let $m\cap \so = \{R_1,R_2,\ldots,R_k\}$.
The lines $R_iN_{\ell}$ for $i=1,2,\ldots,k$ cannot be bisecants of $\so\cup \{N_{\ell}\}$ since they are not tangents to $\so$.
Thus each of these lines meets $\so \cup \{N_{\ell}\}$ in at least three points. Let $B_i=\ell \cap R_iN_{\ell}$, then we have
$|R_iB_i \cap (\so \cup \{N_{\ell}\})|\geq 3$ for $i\in \{1,2,\ldots,k\}$.
We apply Lemma \ref{semioval} part 3 (c) with $T=\ell \cap m$ (note that $k>(a+1)/2\geq 2$).
For $j\in \{2,\ldots,k\}$ we obtain $|R_1B_j \cap (\so \cup \{N_{\ell}\})|=|R_jB_j \cap (\so \cup \{N_{\ell}\})|$,
thus also $|R_1B_j \cap (\so \cup \{N_{\ell}\})| \geq 3$ for $j\in \{2,3,\ldots,k\}$.
We have $N_{\ell}\in R_1B_1$ and hence $N_{\ell}\notin R_1B_j$ for $j\in \{2,3,\ldots,k\}$.
It follows that $R_1B_2 \cup R_1B_3 \cup \ldots R_1B_k \cup m$ contains at least $2(k-1)+k=3k-2$ points of $\so$.
As there is a unique tangent to $\so$ at $R_1$, we must have $a+(q-1)-(3k-2) \geq q-k$. This is a contradiction when $k>(a+1)/2$.
It follows that lines meeting $\so$ in more than $(a+1)/2$ points have to pass through $N_{\ell}$.

Now suppose that $m$ is a $k$-secant of $\so$ with $(a-1)/2 <k<a$ and $N_{\ell}\in m$. Take a point $R\in m\cap \so$.
As $k<a$, there is at least one other line $m'$ through $R$ meeting $\so$ in at least three points.
Let $R'\in (m'\cap \so) \setminus \{R\}$. Lemma \ref{semioval} part 3 (c) with $T=m'\cap \ell$ and $M=m\cap \ell$ yields that
the line joining $R'$ and $m\cap \ell$ meets $\so$ in $|(\so\cup \{N_{\ell}\}) \cap m|=k+1 > (a+1)/2$ points.
Then, according to the previous paragraph, this line also passes through $N_{\ell}$, a contradiction.
It follows that either $k=a$ and hence $N_l \in m$, or $N_l\notin m$ and hence $(a-1)/2<k\leq(a+1)/2$.
\qed

\begin{lemma}
\label{trivi}
Let $\so$ be a semioval of size $q-1+a$ in $\PG(2,q)$.
For each point $R\in \so$ the number of lines through $R$ meeting $\so$ in at least three points is at most $a-2$. \qed
\end{lemma}

\begin{theorem}
\label{2large}
Let $\so$ be a semioval of size $q-1+a$, $a>2$, in $\PG(2,q)$. If $\so$ has two $a$-secants, then one of the following holds.
\begin{enumerate}
	\item $\so$ is the symmetric difference of two lines with one further point removed from both lines.
	\item $\so$ is projectively equivalent to Example \ref{Blokhuis}.
\end{enumerate}
\end{theorem}
\proof
Let $\ell_1$ and $\ell_2$ be two $a$-secants of $\so$ and let $\so'=\so \setminus (\ell_1 \cup \ell_2)$.
Theorem \ref{semioval2} yields $N_1\in\ell_2$ and $N_2\in \ell_1$.
If $\so'=\emptyset$, then $\so \subseteq \ell_1 \cup \ell_2$ and it is easy to see that $\so$ is as in part 1.
If $\so'\neq \emptyset$, then take any point $R\in \so'$. We show that the tangent to $\so$ at $R$ passes through $P:=\ell_1 \cap \ell_2$. As $a>2$, there is a line $r$ through $R$ meeting $\so$ in at least 3 points. According to Theorem \ref{semioval2} part 2, the tangents to $\so$ at the points of $r\cap \so$ pass through a unique point of $\ell_1$, and also through a unique point of $\ell_2$. It follows that these tangents pass through the point $P$.

We show that $\so'$ is contained in the line $\ell_3:=N_1N_2$. Suppose, contrary to our claim, that there is a point $R\in \so' \setminus \ell_3$.
There is a line $r$ through $R$ meeting $\so$ in at least three points. Since $R\notin \ell_3$, $r$ cannot be incident with both $N_1$ and $N_2$.
We may assume $N_2 \notin r$. Let $M=r \cap \ell_1$. Note that $M\notin \so \cup \{N_2,P\}$. Take a point $Q\in \ell_2 \cap \so$. Since the unique tangent to $\so$ at $Q$ is $QN_2$, it follows that $QM$ is a bisecant of $\so$ and it contains a unique point of $\so'$. Denote this point by $R'$.
The tangents to $\so$ at $R$ and $R'$ pass through the same point of $\ell_1$, namely $P$, and hence we can apply
Lemma \ref{semioval} part 3 (d). It follows that $2=|MR' \cap (\so \cup \{N_1\})|=|MR\cap (\so \cup \{N_1\})|\geq 3$.
This contradiction shows $\so' \subset \ell_3$. Lines meeting each of $\ell_1$, $\ell_2$ and $\ell_3$ meet $\so$ in at most two points. Take any point $H\in \so \cap \ell_3$. Since the tangent to $\so$ at $H$ is $PH$, and the other lines through $H$ are not tangents, we obtain $2a=|\ell_1\cap \so|+|\ell_2\cap \so|=q-1$ and hence $a=(q-1)/2$.
The size of $\so$ is $q-1+a=2a+|\so'|$, so $|\so'|=a=(q-1)/2$.
It is easy to show that $\so$ is projectively equivalent to Example \ref{Blokhuis}.
For the complete description of semiovals contained in the sides of a vertexless triangle see the paper of Kiss and Ruff \cite{kissruff}.
\qed
\medskip

A \emph{$(k,n)$-arc} of $\PG(2,q)$ is a set of $k$ points such that each line meets the $k$-set in at most $n$ points.

\begin{theorem}
Let $\so$ be a semioval of size $q+3$ in $\PG(2,q)$, $q$ is a power of the prime $p$. Then $q=5$ and $\so$ is the symmetric difference of two lines with one further point removed from both lines, or $q=9$ and $\so$ is as in Example \ref{Blokhuis}, or $p=3$ and $\so$ is a $(q+3,3)$-arc.
\end{theorem}
\proof
It is easy to see that the points of $\so$ fall into the following two types:
\begin{itemize}
\item points contained in a unique 4-secant and in $q-1$ bisecants,
\item points contained in two trisecants and in $q-2$ bisecants.
\end{itemize}
If $\so$ does not have 4-secants, then the number of trisecants of $\so$ is $(q+3)2/3$, thus $3 \divides q$.
Now suppose that $\so$ has a 4-secant, $\ell$. Theorem \ref{semioval2} with $a=4$ yields that $\so$ does not have trisecants.
The assertion follows from Theorem \ref{2large}.
\qed

\section{\texorpdfstring{Small semiovals when $q$ is even}{Small semiovals when q is even}}

We will use the following theorem by Sz\H{o}nyi and Weiner. This result was proved by the so called resultant method.
We say that a line $\ell$ is an {\em odd-secant} (resp. {\em even-secant}) of $\so$ if $|\ell\cap \so|$ is odd (resp. even).
A \emph{set of even type} is a point set $\cH$ such that each line is an even-secant of $\cH$.

\begin{theorem}[Sz\H{o}nyi and Weiner, \cite{setofeven}]
\label{stabeven}
Assume that the point set $\cH$ in $\PG(2,q)$, $16<q$ even, has $\delta$
odd-secants, where $\delta<(\left\lfloor \sqrt{q}\right\rfloor +1)(q+1-\left\lfloor \sqrt{q}\right\rfloor)$.
Then there exists a unique set $\cH'$ of even type, such that $|\cH \Delta \cH'|=\left\lceil \frac{\delta}{q+1}\right\rceil$.
\end{theorem}

As a corollary of the above result, Sz\H{o}nyi and Weiner gave a lower bound on the size of those point sets of $\PG(2,q)$, $16<q$ even, which do not have tangents but have at least one odd-secant, see \cite{setofeven}.
In this section we prove a similar lower bound on the size of non-oval semiovals.

\begin{lemma}
\label{numodd}
Let $\so$ be a semioval in $\Pi_q$, that is, a projective plane of order $q$.
If $|\so|=q+1+\ep$, then $\so$ has at most $|\so|(1+\ep/3)$ odd-secants.
\end{lemma}
\proof
Take $P\in \so$, then there passes exactly one tangent and there pass at most $\ep$ other odd-secants of $\so$ through $P$.
In this way the non-tangent odd-secants have been counted at least three times.
\qed

\begin{corollary}
\label{evensemi}
If $\so$ is a semioval in $\PG(2,q)$, $16<q$ even, and $|\so|\leq q+3\left\lfloor \sqrt{q}\right\rfloor-11$, then $\so$ is an oval.
\end{corollary}
\proof
If $\delta$ denotes the number of odd-secants of $\so$, then Lemma \ref{numodd} yields:
\[\delta\leq(q+3\left\lfloor \sqrt{q}\right\rfloor-11)(\left\lfloor \sqrt{q}\right\rfloor-3)<(\left\lfloor \sqrt{q}\right\rfloor +1)(q-\left\lfloor \sqrt{q}\right\rfloor+1).\]
By Theorem \ref{stabeven} we can construct a set of even type $\cH$ from $\so$ by modifying (add to $\so$ or delete from $\so$)
$\left\lceil \frac{\delta}{q+1}\right\rceil \leq \left\lfloor \sqrt{q}\right\rfloor+1$ points of $\PG(2,q)$.

If $P\in \so$ is a modified (and hence deleted) point, then the number of lines through $P$ which are not tangents to $\so$ and do not
contain modified points is at least $q-\left(\left\lceil \frac{\delta}{q+1}\right\rceil-1\right)$.
These lines are even-secants of $\cH$ and hence they are non-tangent odd-secants of $\so$.
It follows that the size of $\so$ is at least $1+2(q-\left\lfloor \sqrt{q}\right\rfloor)$, a contradiction.

Thus each of the modified points has been added.
Suppose $|\so|>q+1$. As there is a tangent to $\so$ at each point of $\so$, we have $2 \leq \left\lceil \frac{\delta}{q+1}\right\rceil$.
Let $A$ and $B$ be two modified (and hence added) points.
If the line $AB$ contains another added point $C$, then through one of the points $A$, $B$, $C$ there pass at most $(|\so|-1)/3+1$ tangents to $\so$. If $AB$ does not contain further added points, then $AB$ cannot be a tangent to $\so$ and hence through one of the points
$A$, $B$ there pass at most $|\so|/2$ tangents to $\so$. Let $A$ be an added point through which there pass at most $|\so|/2$ tangents to $\so$ and denote the number of these tangents by $\tau$. Through $A$ there pass at least $q+1-\tau-\left(\left\lceil \frac{\delta}{q+1}\right\rceil-1\right)$ lines meeting $\so$ in at least two points. Thus from $\tau \leq |\so|/2$ and from the assumption on the size of $\so$ we get
\[q+3\left\lfloor \sqrt{q}\right\rfloor-11\geq \tau + 2(q+1-\tau-\left\lfloor \sqrt{q}\right\rfloor)
\geq 2(q-\left\lfloor \sqrt{q}\right\rfloor+1)-(q+3\left\lfloor \sqrt{q}\right\rfloor-12)/2.\]
After rearranging we obtain $0 \geq q - 13\left\lfloor \sqrt{q}\right\rfloor + 38$, which is a contradiction.
It follows that $|\so|\leq q+1$, but also $|\so|\geq q+1$ and $\so$ is an oval in the case of equality.
\qed

\section{\texorpdfstring{Point sets with few odd-secants in $\PG(2,q)$, $q$ odd}{Point sets with few odd-secants in PG(2,q), q odd}}

Some combinatorial results of this section hold in every finite projective plane.
As before, by $\Pi_q$ we denote an arbitrary projective plane of order $q$.

\begin{definition}
\label{type}
Fix a point set $\so \subseteq \Pi_q$.
For a positive integer $i$ and a point $P\in \so$ we denote by $t_i(P)$ the number of $i$-secants of $\so$ through $P$.
The weight of $P$, in notation $w(P)$, is defined as follows.
\[w(P):=\sum_{i \text{ odd}}t_i(P)/i.\]
For a subset $\cP \subseteq \so$, let $w(\cP)=\sum_{P\in \cP} w(P)$.
Suppose that $w(P)$ is known for $P \in \{P_1,P_2,\ldots,P_m\}\subseteq \so \cap \ell$, where $\ell$ is a line
meeting $\so$ in at least $m$ points. Then the type of $\ell$ is
\[[w(P_1),w(P_2),\ldots,w(P_m)].\]
Suppose that the value of $t_i(P)$ is known for a point $P\in \so$ and for $1 \leq i\leq q+1$.
Let $\{a_1,a_2,\ldots,a_k\}=\{i \colon t_i(P)\neq 0\}$, then the type of $P$ is
\[[a_{1\,t_{a_1}(P)},a_{2\,t_{a_2}(P)},\ldots,a_{k\,t_{a_k}(P)}].\]
\end{definition}

\begin{example}[Balister et al. \cite{symmdiff}]
Let $\so=\cC \cup \{P\}$, where $\cC$ is a conic of $\PG(2,q)$, $q$ odd, and $P\notin \cC$ is an external point of $\cC$, that is, a point contained in two tangents to $\cC$. Then the type of $P$ is $[1_{(q-1)/2},2_2,3_{(q-1)/2}]$ and $w(P)=(q-1)/2+(q-1)/6$. If $T_1$ and $T_2$ are the points of $\cC$ contained in the tangents to $\cC$ at $P$, then the type of $T_i$ is $[2_{q+1}]$ and $w(T_i)=0$ for $i=1,2$. Each point of $\cC \setminus \{T_1,T_2\}$ has type $[1_1,2_{q-1},3_1]$ and weight $4/3$. The number of odd-secants of $\so$ is $2q-2$.
\end{example}

\begin{theorem}[Balister et al. {\cite[Theorem 6]{symmdiff}}]
The minimal number of odd-secants of a $(q+2)$-set in $\PG(2,q)$, $q$ odd, is $2q-2$ when $q\leq 13$.
For $q\geq 7$, it is at least $3(q+1)/2$.
\end{theorem}

\begin{conjecture}[Balister et al. {\cite[Conjecture 11]{symmdiff}}]
\label{2q-2}
The minimal number of odd-secants of a $(q+2)$-set in $\PG(2,q)$, $q$ odd, is $2q-2$.
\end{conjecture}
The following propositions are straightforward.
\begin{proposition}
The number of odd-secants of $\so$ is $w(\so)=\sum_{P\in \so}w(P)$.
\qed
\end{proposition}

\begin{proposition}
\label{start}
Let $\so$ be a $(q+2)$-set in $\Pi_q$ and let $P$ be a point of $\so$.
The smallest possible weights of $P$ are as follows:
\begin{itemize}
	\item $w(P)=0$ if and only if the type of $P$ is $[2_{q+1}]$,
	\item $w(P)=4/3$ if and only if the type of $P$ is $[1_1,2_{q-1},3_1]$,
	\item $w(P)=2$ if and only if the type of $P$ is $[1_2,2_{q-2},4_1]$,
	\item $w(P)=8/3$ if and only if the type of $P$ is $[1_2,2_{q-3},3_2]$,
	\item $w(P)=16/5$ if and only if the type of $P$ is $[1_3,2_{q-2},5_1]$,
	\item $w(P)=10/3$ if and only if the type of $P$ is $[1_3,2_{q-3},3_1,4_1]$. \qed
\end{itemize}
\end{proposition}

\begin{proposition}
\label{start2}
Let $\so$ be a point set of size $q+2$ in $\Pi_q$ and let $P$ be a point of $\so$.
\begin{enumerate}
	\item If $P$ is contained in a $k$-secant, then $w(P)\geq k-2$,
	\item if $P$ is contained in at least $k$ trisecants, then $w(P)\geq \frac{4}{3}k$.
\end{enumerate}
\end{proposition}
\proof
In part 1, the number of tangents to $\so$ at $P$ is at least $q-(q+2-k)=k-2$.
In part 2, $P$ is incident with at least $q+1-k-(q+2-(2k+1))=k$ tangents to $\so$, thus $w(P)\geq k/3+k$.
\qed

\begin{theorem}[Bichara and Korchm\'{a}ros {\cite[Theorem 1]{34bakg}}]
\label{bk}
Let $\so$ be a point set of size $q+2$ in $\PG(2,q)$. If $q$ is odd, then $\so$ contains at most two points with weight 0, that is, points of type $[2_{q+1}]$.
\end{theorem}

\begin{lemma}
\label{nuk}
Let $\so$ be a point set of size $q+k$ in $\PG(2,q)$ for some $k\geq 3$.
Suppose that $\ell_1$ is a $k$-secant of $\so$ meeting $\so$ only in points of type $[2_{q},k_1]$.
Then the $k$-secants of $\so$ containing a point of type $[2_{q},k_1]$ are concurrent.
\end{lemma}
\proof
Let $\ell_2, \ell_3$ be two $k$-secants of $\so$ with the given property and let $R_i \in \ell_i\cap \so$ be a point of type $[2_q,k_1]$ for $i=2,3$. It is easy to see that $\cB:=\ell \Delta \so$ is a blocking set of R\'edei type and $R_2$, $R_3$ are not incident with any bisecant of $\cB$. It follows from Theorem \ref{main0} part 2 that $\ell_2 \cap \ell_3 \in \ell_1$.
\qed

\begin{definition}
A \emph{$(q+t,t)$-arc of type $(0,2,t)$} is a point set $\cT$ of size $(q+t)$ in $\PG(2,q)$ such that each line meets $\cT$ in 0,2 or $t$ points.
In honor of Korchm\'aros and Mazzocca such point sets are also called \emph{KM-arcs} in the literature.
\end{definition}

Let $\cT$ be a $(q+t,t)$-arc of type $(0,2,t)$.
It is easy to see that for $t>2$ there is a unique $t$-secant through each point of $\cT$. It can be proved that $2\leq t <q$ implies $q$ even, see \cite{KMarcs} by Korchm\'aros and Mazzocca. As the points of $\cT$ are of type $[2_q,t_1]$, the following theorem by G\'acs and Weiner also follows from
Lemma \ref{nuk}. For recent results on KM-arcs we refer the reader to \cite{KM-arcs}.

\begin{theorem}[G\'acs and Weiner {\cite[Theorem 2.5]{gw}}]
Let $\cT$ be a $(q+t,t)$-arc of type $(0,2,t)$ in $\PG(2,q)$. If $t>2$, then the $t$-secants of $\cT$ pass through a unique point. \qed
\end{theorem}

The proof of our next result is based on the counting technique of Segre. A \emph{dual arc} is a set of lines such that no three of them are concurrent.

\begin{theorem}
\label{gbk}
Let $\so$ be a point set of size $q+k$ in $\PG(2,q)$, $q$ odd.
\begin{enumerate}
	\item If $k=1$, then the tangents to $\so$ at points of type $[1_1,2_q]$ form a dual arc.
	\item If $k=2$, then there are at most two points of type $[2_{q+1}]$.
	\item If $k\geq 3$, then the $k$-secants of $\so$ containing a point of type $[2_{q},k_1]$ form a dual arc.
\end{enumerate}
\end{theorem}
\proof
Suppose the contrary. If $k=1$, then let $A$, $B$ and $C$ be points of type $[1_1,2_q]$ such
that the tangents through these points pass through a common point $D$.
If $k=2$, then let $A$, $B$ and $C$ be three points of type $[2_{q+1}]$ and take a point $D\notin (\so \cup AB \cup BC \cup CA)$.
If $k\geq 3$, then let $A$, $B$ and $C$ be points of type $[2_q,k_1]$ such that the $k$-secants through these points pass through a common point
$D\notin (AB \cup BC \cup CA)$.
In all cases, $A$, $B$, $C$ and $D$ are in general position, thus we may assume $A=(\infty)$, $B=(0,0)$, $C=(0)$ and $D=(1,1)$.
Let $\so'=\so \setminus \{A,B,C\}$. Note that $AB$, $BC$ and $CA$ are bisecants of $\so$ and $CA$ is the line at infinity, thus $\so'$
is a set of $q+k-3$ affine points, say $\so'=\{(a_i,b_i)\}_{i=1}^{q+k-3}$.
For $i\in \{1,2,\ldots,q+k-3\}$ we have the following.
\begin{itemize}
	\item The line joining $(a_i,b_i)$ and $A$ meets $BC$ in $(a_i,0)$,
	\item the line joining $(a_i,b_i)$ and $B$ meets $AC$ in $(b_i/a_i)$,
	\item the line joining $(a_i,b_i)$ and $C$ meets $AB$ in $(0,b_i)$.
\end{itemize}
The lines $AD$, $BD$ and $CD$ meet $\so'$ in $k-1$ points.
The lines $AP$ for $P\in \so' \setminus AD$ meet $\so'$ in a unique point.
Since the first coordinate of the points of $AD\cap \so'$ is 1, it follows that $\{a_i\}_{i=1}^{q+k-3}$ is a multiset containing each element of $\GF(q)\setminus \{0,1\}$ once, and containing $1$ $k-1$ times. Thus $\prod_{i=1}^{q+k-3}a_i=-1$.
Similarly, the lines through $B$ yield $\prod_{i=1}^{q+k-3}b_i/a_i=-1$, and the lines through $C$ yield $\prod_{i=1}^{q+k-3}b_i=-1$.
It follows that
\[1=(-1)(-1)=\left(\prod_{i=1}^{q+k-3}a_i\right)\left(\prod_{i=1}^{q+k-3}\frac{b_i}{a_i}\right)=\prod_{i=1}^{q+k-3}b_i=-1,\]
a contradiction for odd $q$.
\qed
\medskip

The following immediate consequence of Theorem \ref{gbk} and Lemma \ref{nuk} will be used frequently.

\begin{corollary}
\label{all}
Let $\so$ be a point set of size $q+k$, $k\geq 3$, in $\PG(2,q)$. If there exist three $k$-secants of $\so$, $\ell_1$, $\ell_2$ and $\ell_3$,
such that the points of $\ell_1 \cap \so$ are of type $[2_q,k_1]$ and both $\ell_2 \cap \so$ and $\ell_3 \cap \so$ contain at least one point of type $[2_q,k_1]$, then $q$ is even.
\end{corollary}
\proof
Lemma \ref{nuk} yields $\ell_2 \cap \ell_3 \in \ell_1$, but then Theorem \ref{gbk} implies $q$ even.
\qed

For the definition of a nucleus $N_i$ of a line $\ell_i$ see Proposition \ref{nucleus}.

\begin{lemma}
\label{elsocor}
Let $\so$ be a set of $q-1+a$ points, $a\geq 3$, in $\PG(2,q)$, where $q$ is a power of the prime $p$. Suppose that $\ell_1$ and $\ell_2$ are $a$-secants of $\so$ such that there is a unique tangent to $\so$ at each point of $\so \cap \ell_i$, for $i=1,2$.
\begin{enumerate}
	\item Either $N_1\in \ell_2$ and $N_2 \in \ell_1$, or
	\item $N_1=N_2$, $p \divides a$ and for each $R\in\so$ if there is a unique tangent $r$ to $\so$ at $R$, then $r$ passes through the common nucleus.
	\item Let $\ell_3$ be another $a$-secant of $\so$ such that there is a unique tangent to $\so$ at each point of $\so \cap \ell_3$.
	If $q$ or $a$ is odd, then $\ell_3=N_1N_2$, thus in this case $\ell_3$ is uniquely determined.
\end{enumerate}
\end{lemma}
\proof
If $\ell_1 \cap \ell_2 \in \so$, then $|\so|\geq 2a+q-3$, which cannot be since $a\geq 3$. First assume $N_1 \neq N_2$ and suppose to the contrary $N_2 \notin \ell_1$. Then $\cB:=\{N_1\} \cup (\ell_1 \Delta \so)$ is a blocking set of R\'edei type.
There is a unique bisecant of $\cB$ at each point of $\so \cap \ell_2$ (the tangent to $\so$). This is a contradiction since these bisecants should pass through the same point of $\ell_1$ (apply Theorem \ref{main1} part 2 with $T=\ell_1\cap \ell_2$).

If $N_1=N_2=:N$, then we define $\cB$ in the same way.
Then there is no bisecant of $\cB$ through the points of $\cB \cap \ell_2$. Theorem \ref{main0} yields $p \divides a$.
Take a point $R\in \so \setminus (\ell_1 \cup \ell_2)$ incident with a unique tangent $r$ to $\so$. If $N \notin r$, then $r$ is the unique bisecant  of $\cB$ through $R$, a contradiction because of Theorem \ref{main1} part 1.

Suppose that $\ell_3$ is an $a$-secant with properties as in part 3. Then either $\ell_3=N_1N_2$ and $N_3=\ell_1 \cap \ell_2$, or
$N_3=N_1=N_2=:N$ and $p\divides a$. In the latter case Corollary \ref{all} applied to $\so \cup \{N\}$ and to the lines
$\ell_1$, $\ell_2$ and $\ell_3$ yields $p=2$.
\qed

\begin{lemma}
\label{masodikcor}
Let $\so$ be a set of $q+2$ points in $\PG(2,q)$, $q$ is a power of the odd prime $p$, and suppose that $\ell$ is a trisecant of $\so$ of type $[4/3,4/3,4/3]$.
\begin{enumerate}
	\item If $p=3$, then the tangents at the points of $\so$ with weight 4/3 pass through $N_{\ell}$.
	There is at most one other trisecant of $\so$ of type $[4/3]$.
	\item If $p\neq 3$, then the trisecants of type $[4/3,4/3]$ pass through $N_{\ell}$.
	Suppose that there is another trisecant $\ell_1$ of type $[4/3,4/3,4/3]$. Then there is at most one other trisecant of type
	$[4/3,4/3]$, which is $N_{\ell}N_1$.
	If $N_{\ell}N_1$ is a trisecant of type $[4/3,4/3]$, then the tangents at the points of $N_{\ell}N_1$ with weight 4/3
	pass through $\ell\cap \ell_1$.
\end{enumerate}
\end{lemma}
\proof
Let $\cB$ denote the R\'edei type blocking set $(\ell \Delta \so) \cup \{N_{\ell}\}$.

First we prove part 1. Take $A\in \so \setminus \ell$ such that $w(A)=4/3$ and denote the tangent to $\so$ at $A$ by $a$.
If $N_{\ell}\notin a$, then there is a unique bisecant of $\cB$ through $A$, thus Theorem \ref{main1} yields $p\neq 3$, a contradiction.
Denote the trisecant through $A$ by $\ell_1$.
If there were a trisecant $\ell_2$ of type $[4/3]$ different from $\ell$ and $\ell_1$, then
Corollary \ref{all} applied to $\so \cup \{N_{\ell}\}$ and to the lines
$\ell$, $\ell_1$ and $\ell_2$ would yield $q$ even, a contradiction.

Now we prove part 2. First suppose to the contrary that there is a trisecant $\ell_2$ of type $[4/3,4/3]$ with $N_{\ell}\notin \ell_2$.
Let $A,B\in \ell_2 \cap \so$ such that $w(A)=w(B)=4/3$. Denote the tangents to $\so$ at these two points by $a$ and $b$, respectively.
We have $N_{\ell} \notin a$ and $N_{\ell}\notin b$, since otherwise we would get points not incident with any bisecant of $\cB$,
a contradiction as $p\neq 3$ (cf. Theorem \ref{main0}). It follows that $N_{\ell}A$ and $N_{\ell}B$ are 4-secants of $\cB$. Let $M=N_{\ell}A \cap \ell$.
Then Theorem \ref{main1} part 2 (with $T=\ell \cap \ell_2$) yields that $MB$ is also a 4-secant of $\cB$ and hence a trisecant of $\so$ (we have $N_{\ell}\notin MB$). A contradiction, since $MB\neq \ell_2$. It follows that $N_{\ell}\in \ell_2$.

Let $\ell_1$ be trisecant of $\so$ of type $[4/3,4/3,4/3]$ and let $\ell_2$, $A$, $B$, $a$ and $b$ be defined as in the previous paragraph.
It follows from Lemma \ref{elsocor} that $N_{\ell}\in \ell_1$ and $N_1\in \ell$.
It also follows from the previous paragraph that $N_1 \in \ell_2$ and $N_{\ell}\in \ell_2$, thus $\ell_2=N_1N_{\ell}$.
Theorem \ref{main1} applied to $\cB$ and to $(\ell_1 \Delta \so)\cup \{N_1\}$ yields that $a$ and $b$ pass through a unique point of $\ell$ and through a unique point of $\ell_1$, thus they pass through $\ell \cap \ell_1$.
\qed
\medskip

Let $\so$ be a set of $q+2$ points of $\PG(2,q)$, $q$ odd. Since $q+2$ is odd, each point $P\notin \so$ is incident with an odd-secant of $\so$. It follows that the odd-secants of $\so$ cover the points of $\PG(2,q)$ except for the points of $\so$ with weight zero. For partial covers of $\PG(2,q)$ we refer the reader to \cite[Proposition 1.5]{partial}. The lower bound on the size of an affine blocking set \cite{BS, J} yields the following result. Its proof can be found in \cite{kakeya} at the top of page 211, as part of a more complex argument. For a proof in the dual setting see \cite[Lemma 10]{symmdiff}.

\begin{lemma}[Blokhuis and Mazzocca {\cite{kakeya}}]
\label{kaki}
Let $\so$ be a set of $q+2$ points of $\PG(2,q)$, $q$ odd. If $\so$ has $d\in\{1,2\}$ points with weight zero, then the number of odd-secants of
$\so$ is at least $2q-d$.
\end{lemma}

\begin{theorem}
Let $\so$ be a point set of size $q+2$ in $\PG(2,q)$, $13 < q$ odd.
Then the number of odd-secants of $\so$ is at least $\left\lceil\frac{8}{5}q + \frac{12}{5}\right\rceil$.
\end{theorem}
\proof
Let $d$ denote the number of points of $\so$ with weight zero. Theorem \ref{bk} of Bichara and Korchm\'aros yields $d\leq 2$.
If $d\in \{1,2\}$, then Lemma \ref{kaki} yields $w(\so)\geq 2q-2$, which is at least $\left\lceil\frac{8}{5}q + \frac{12}{5}\right\rceil$ when $q\geq 11$.
From now on we assume $d=0$. Consider the following subsets of $\so$:
\[\cB:=\{P\in \so \colon \text{$P$ is contained in a trisecant of type [4/3,4/3,4/3]}\},\]
\[\cC:=\{P\in \so \colon w(P)\neq 4/3,\, \text{$P$ is contained in a trisecant of type [4/3]}\}.\]
Denote the size of $\cC$ by $m$ and let $\cC=\{P_1,P_2,\ldots,P_m\}$.
For $i=1,2,\ldots,m$, let
\[V_i= \{Q\in \so  \colon w(Q)=4/3 \text{ and $QP_i$ is a trisecant}\} \cup \{P_i\}.\]
Also, let $D_1:=V_1$ and $D_i:=V_i \setminus (\cup_{j=1}^{i-1}V_j)$ for $i\in \{2,3,\ldots,m\}$.
Of course the sets $D_1,D_2,\ldots,D_m$ are disjoint and $P_i\in D_i\subseteq V_i$.
The point set $\cD:=\cup_{i=1}^m D_i$ contains each point of $\so \setminus \cB$ with weight 4/3. Note that each point of $D_i$
has weight 4/3, except $P_i$. We introduce the following notion. For a point set $\cU\subseteq \so$ let $\alpha(\cU)$ denote the average weight of the points in $\cU$, that is, $\alpha(\cU)=w(\cU)/|\cU|$. First we prove $\alpha(D_i)\geq 8/5$ for $i=1,2,\ldots,m$.
If $t_3(P_i)=k$ (cf. Definition \ref{type}), then
\begin{equation}
\label{2k+1}
|D_i|\leq|V_i|\leq 2k+1.
\end{equation}
If $k=1$, then Proposition \ref{start} yields $w(P_i)\geq 10/3$ (since $w(P_i)\neq 4/3$), hence in this case we have
\begin{equation}
\label{av1}
\alpha(D_i)\geq \frac{10/3+(|D_i|-1)4/3}{|D_i|}=4/3+\frac{2}{|D_i|}\geq 2.
\end{equation}
If $k\geq 2$, then Proposition \ref{start2} yields $w(P_i)\geq 4k/3$, thus
\begin{equation}
\label{av2}
\alpha(D_i)\geq \frac{4k/3+(|D_i|-1)4/3}{|D_i|}=4/3+\frac{(k-1)4/3}{|D_i|}\geq 2-\frac{2}{2k+1}\geq 8/5.
\end{equation}
We define a further subset of $\so$, $\cE:=\so \setminus (\cB \cup \cD)$.
Note that $w(\cD)\geq |\cD|\frac85$ and $w(\cE)\geq |\cE|2$, since each point of $\cE$ has weight at least 2 (see Porposition \ref{start}).
The point sets $\cB$, $\cD$ and $\cE$ form a partition of $\so$, thus $w(\so)=w(\cB)+w(\cD)+w(\cE)$.
We distinguish three main cases.
\begin{enumerate}
	\item There is no trisecant of $\so$ of type $[4/3,4/3,4/3]$. Then we obtain $w(\so)\geq (q+2)\frac85$.
	\item There is at least one trisecant of $\so$ of type $[4/3,4/3,4/3]$ and $p\neq 3$.
		Denote the number of trisecants of $\so$ of type $[4/3,4/3,4/3]$ by $s$. Lemma \ref{masodikcor} yields $s\leq 3$.
		If $s=1$, then $w(s)\geq 3\frac43+(q-1)\frac85=q\frac85+\frac{12}{5}$.
		If $s=2$, then according to Lemma \ref{masodikcor} there is at most one other trisecant of type $[4/3,4/3]$.
		Thus in \eqref{2k+1} we have $|D_i|\leq |V_i|\leq k+2$, where $k=t_3(P_i)$.
		If $k=1$, then similarly to \eqref{av1} we obtain $\alpha(D_i)\geq 2$.
		If $k\geq 2$, then similarly to \eqref{av2} we obtain $\alpha(D_i)\geq \frac53$. It follows that
		$w(\so)\geq 6\frac43+(q-4)\frac53=q\frac53+\frac43$.
		If $s=3$, then according to Lemma \ref{masodikcor} there is no other trisecant of type $[4/3,4/3]$.
		Thus in \eqref{2k+1} we have $|D_i|\leq |V_i|\leq k+1$. If $k=1$, then similarly to \eqref{av1} we obtain
		$\alpha(D_i)\geq \frac73$, if $k\geq 2$, then similarly to \eqref{av2} we obtain $\alpha(D_i)\geq \frac{16}{9}$. It follows that
		$w(\so)\geq 9\frac43+(q-7)\frac{16}{9}=q\frac{16}{9}-\frac49$.
	\item There is at least one trisecant $\ell$ of $\so$ of type $[4/3,4/3,4/3]$ and $p=3$. It follows from Lemma \ref{masodikcor} that
		the number $g$ of further trisecants of type $[4/3]$ is at most one. First suppose $g=0$. As $\cD$ is empty, we obtain $w(\so)\geq 3\frac43+(q-1)2\geq 2q+2$. If $g=1$, then let $r\neq \ell$ be the other trisecant of $\so$ of type [4/3]. Let $t\in \{1,2,3\}$ be
		the number of points with weight 4/3 in $r\cap \so$. It follows that
		$w(\so)\geq (3+t)\frac43+(3-t)\frac83+(q-4)2\geq 6\frac43+(q-4)2=2q$.
\qed
\end{enumerate}

For a line set $\cL$ of $\AG(2,q)$, $q$ odd, denote by $\tilde{w}(\cL)$ the set of affine points contained in an odd number of lines of $\cL$.
\cite[Theorem 3.2]{PV} by Vandendriessche classifies those line sets $\cL$ of $\AG(2,q)$ for which $|\cL|+\tilde{w}(\cL)\leq 2q$, except for
one open case ({\cite[Open Problem 3.3]{PV}}), which we recall here. For applications in coding theory we refer the reader to the Introduction of the paper of Vandendriessche and the references there.

\begin{example}[Vandendriessche {\cite[Example 3.1 (i)]{PV}}]
\label{Peter}
$\cL$ is a set of $q+k$ lines in $\AG(2,q)$, $q$ odd, with the following properties. There is an
$m$-set $\so \subset \ell_{\infty}$ with $4\leq m \leq q-1$ and an odd positive integer $k$ such that
exactly $k$ lines of $\cL$ pass through each point of $\so$ and $\tilde{w}(\cL)=q-k$.
\end{example}

\begin{proposition}
Example \ref{Peter} cannot exist.
\end{proposition}
\proof
The dual of the line set $\cL$ in Example \ref{Peter} is a point set $\cB$ of size $q+k$ in $\PG(2,q)$, such that there is a point $O\notin \cB$ (corresponding to $\ell_{\infty}$), with the properties that through $O$ there pass $m$ $k$-secants of $\cB$, $\ell_1, \ell_2,\ldots,\ell_m$, and the number of odd-secants of $\cB$ not containing $O$ is $q-k$ ($q$, $m$ and $k$ are as in Example \ref{Peter}).

As $q+k$ is even and $k$ is odd, it follows for $i\in \{1,2,\ldots,m\}$ and for any $R\in \ell_i \setminus (\cB \cup \{O\})$ that through $R$ there passes at least one odd-secant of $\cB$, which is different from $\ell_i$. As the number of odd-secants of $\cB$ not containing $O$ is $q-k$, and $|\ell_i \setminus (\cB \cup \{O\})|=q-k$, it follows that there is a unique odd-secant of $\cB$ through each point of $\cB \cap \ell_i$, namely $\ell_i$.
But $|\cB \setminus \ell_i|=q$, thus lines not containing $O$ and meeting $\ell_i$ in $\cB$ are bisecants of $\cB$ (otherwise we would get tangents to $\cB$ not containing $O$ at some point of $\ell_i\cap \cB$). Then for $i\in \{1,2,\ldots,m\}$ the points of $\cB\cap \ell_i$ are of type $[2_q,k_1]$. As $m\geq 3$ and the lines $\ell_1,\ldots, \ell_m$ are concurrent, Theorem \ref{gbk} yields a contradiction for odd $q$.
\qed

\begin{remark}
Together with other ideas, our method yields lower bounds on number of odd-secants of $(q+3)$-sets and $(q+4)$-sets as well.
We will present these results elsewhere.
\end{remark}

{\bf Acknowledgement.} The author is grateful to the referees for their useful comments, in particular for the insight that a previous version of Theorem \ref{2large} can be improved.

\begin{flushleft}
Bence Csajb\'ok \\
Dipartimento di Tecnica e Gestione dei Sistemi Industriali,\\
Universit\`a di Padova,\\
Stradella S. Nicola, 3, I-36100 Vicenza, Italy\\
and \\
MTA--ELTE Geometric and Algebraic Combinatorics Research Group, \\
E\"otv\"os Lor\'and University, \\
1117 Budapest, P\'azm\'any P\'eter S\'et\'any 1/C, Hungary, \\
e-mail: {\sf csajbok.bence@gmail.com}
\end{flushleft}


\begin{thebibliography}{10}

\bibitem{symmdiff}
P.~Balister, B.~Bollob\'{a}s, Z.~F\"{u}redi and J.~Thompson,
{\it Minimal Symmetric Differences of Lines in Projective Planes},
J. Combin. Des. 22(10) (2014), 435--451.

\bibitem{Ball}
S.~Ball, {\it The number of directions determined by a function over a finite field}, J. Combin. Theory Ser. A
104 (2003), 341--350.

\bibitem{BBBSSz}
S.~Ball, A.~Blokhuis, A.E.~Brouwer, L.~Storme and T.~Sz\H{o}nyi,
{\it On the number of slopes of the graph of a function definied over a finite field}, J. Combin. Theory Ser. A
86 (1999), 187--196.

\bibitem{db}
D.~Bartoli, {\it On the Structure of Semiovals of Small Size},
J. Combin. Des. 22(12) (2014), 525–-536.

\bibitem{34bakg}
A.~Bichara and G.~Korchm\'aros, {\it Note on $(q+2)$-sets in a Galois plane of order $q$},
Ann. Discrete Math. 14 (1980), 117--121.

\bibitem{seminuclear}
A.~Blokhuis, {\it Characterization of seminuclear sets in a finite projective plane},
J. Geom. 40 (1991), 15--19.

\bibitem{BB}
A.~Blokhuis and A.E.~Brouwer, {\it Blocking sets in Desarguesian projective planes}, Bull. London Math. Soc. 18 (1986), 132–-134.

\bibitem{partial}
A.~Blokhuis, A.E.~Brouwer and T. Sz\H{o}nyi, {\it Covering all points except one}, J. Algebraic Combin. 32 (2010), 59--66.

\bibitem{originsemi}
A.~Blokhuis and A.A.~Bruen, {\it The minimal number of lines intersected by a set of $q+2$ points, blocking sets and intersecting circles}, J. Combin. Theory Ser. A 50 (1989), 308–-315.

\bibitem{kakeya}
A.~Blokhuis and F.~Mazzocca, {\it The finite field Kakeya problem}, Building Bridges 205--218, Bolyai Soc. Math. Stud. 19, Springer, Berlin, 2008.

\bibitem{BS}
A.E.~Brouwer and A. Schrijver, {\it The blocking number of an affine space}, J. Combin. Theory Ser. A 24 (1978), 251--253.

\bibitem{semiar2}
B.~Csajb\'{o}k, T.~H\'{e}ger and Gy.~Kiss, {\it Semiarcs with a long secant in $\mathrm{PG}(2,q)$}, Innov. Incidence Geom. 14 (2015), 1--26.

\bibitem{KM-arcs}
M.~De Boeck and G. Van de Voorde, {\it A linear set view on KM-arcs}, to appear in J. Algebraic Combin., DOI 10.1007/s10801-015-0661-7

\bibitem{perpol}
R.J.~Evans, J.~Greene, H.~Niederreiter,
{\it Linearized polynomials and permutation polynomials of finite fields}, Michigan Math. J. 39 (1992), 405--413.

\bibitem{gacssemi}
A.~G\'{a}cs, {\it On regular semiovals in $\mathrm{PG}(2,q)$}, J. Algebraic Combin. 23 (2006), 71--77.

\bibitem{gw}
A.~G\'{a}cs and Zs.~Weiner, {\it On $(q+t,t)$-arcs of type $(0,2,t)$}, Des. Codes Cryptogr. 29 (2003), 131--139.

\bibitem{jwph}
J.W.P.~Hirschfeld, {\it Projective Geometries over Finite Fields,} 2$^{nd}$ ed., Clarendon Press, Oxford, 1998.

\bibitem{J}
R.~Jamison, {\it Covering finite fields with cosets of subspaces}, J. Combin. Theory Ser. A 22 (1977), 253--266.

\bibitem{survey}
Gy.~Kiss, {\it A survey on semiovals}, Contrib. Discrete Math. 3 (2008), 81--95.

\bibitem{smalls}
Gy.~Kiss, S.~Marcugini and F.~Pambianco, {\it On the spectrum of the sizes of semiovals in $\mathrm{PG}(2,q)$, $q$ odd},
Discrete Math. 310 (2010), 3188--3193.

\bibitem{kissruff}
Gy.~Kiss and J.~Ruff, {\it Notes on Small Semiovals}, Annales Univ. Sci. Budapest 47 (2004), 143--151.

\bibitem{KMarcs}
G.~Korchm\'aros and F.~Mazzocca, {\it On $(q+t)$-arcs of type $(0,2,t)$ in a desarguesian plane of order $q$},
Math. Proc. Cambridge Philos. Soc. 108 (1990), 445--459.

\bibitem{lisonek}
P.~Lisonek, Computer-assisted Studies in Algebraic Combinatorics, Ph.D. Thesis, RISC, J. Kepler University Linz, 1994.

\bibitem{linconj}
P.~Sziklai, {\it On small blocking sets and their linearity}, J. Combin. Theory Ser. A 115 (2008), 1167--1182.

\bibitem{Szonyiblok}
T.~Sz\H{o}nyi, {\it Blocking Sets in Desarguesian Affine and Projective Planes}, Finite Fields Appl. 3 (1997), 187--202.

\bibitem{dirszonyi}
T.~Sz\H{o}nyi, {\it On the Number of Directions Determined by a Set of Points in an Affine Galois Plane},
J. Combin. Theory Ser. A 74 (1996), 141--146.

\bibitem{setofeven}
T.~Sz\H{o}nyi and Zs.~Weiner, {\it On the stability of the sets of even type},
Adv. Math. 267 (2014), 381-–394.

\bibitem{PV}
P.~Vandendriessche, {\it On small line sets with few odd-points},
Des. Codes Cryptogr. 75 (2015), 453--463.

\end{thebibliography}
\end{document}